\documentclass{amsart}

\usepackage{amsmath,amssymb}
\usepackage{hyperref,url}
\usepackage{tikz}
\usepackage{colortbl}
\usepackage{caption}

\usetikzlibrary{matrix}
\usetikzlibrary{positioning}
\usetikzlibrary{calc}

\begin{document}

\title{Generating modular lattices of up to 30 elements}

\author{Jukka Kohonen}
\address{Department of Computer Science, Aalto University, Espoo, Finland.
  \emph{Present address:} Department of Computer Science, University of Helsinki, Helsinki, Finland}
\email{jukka.kohonen@iki.fi}

\begin{abstract}
  An algorithm is presented for generating finite modular,
  semimodular, graded, and geometric lattices up to isomorphism.
  Isomorphic copies are avoided using a combination of the
  general-purpose graph-isomorphism tool \texttt{nauty} and some
  optimizations that handle simple cases directly.  For modular and
  semimodular lattices, the algorithm prunes the search tree much
  earlier than the method of Jipsen and Lawless, leading to a speedup
  of several orders of magnitude.  With this new algorithm modular
  lattices are counted up to 30 elements, semimodular lattices up to
  25 elements, graded lattices up to 21 elements, and geometric
  lattices up to 34 elements.  Some statistics are also provided on
  the typical shape of small lattices of these types.
\end{abstract}
\keywords{Modular lattices, semimodular lattices, graded lattices, geometric lattices, counting algorithm}

\maketitle

\section{Introduction}
\label{sec:intro}

Algorithms that generate and count unlabeled lattices follow generally
the same pattern: start from a small initial lattice, recursively add
new elements, and take care to keep only one representative of each
isomorphism class.  With variations of this scheme, unlabeled lattices
have been counted up to $18$~elements by Heitzig and
Reinhold~\cite{heitzig2002}, $19$~elements by Jipsen and
Lawless~\cite{jipsen2015}, and $20$~elements by Gebhardt and
Tawn~\cite{gebhardt2016}.  All of these enumerations took hundreds of
days of processor time.  Special lattices may be generated faster, or
to a larger size: Jipsen and Lawless counted \emph{modular} lattices
up to $24$ elements and \emph{semimodular} lattices up to $22$
elements~\cite{jipsen2015}.  Empirically the running time of their
method grows faster than~${6}^n$, where $n$ is lattice size (number of
elements).

This paper describes an improved algorithm for generating graded
lattices and certain subfamilies.  In easy cases it handles
isomorphisms quickly, avoiding a costly call to \texttt{nauty}.  But
more importantly, with (semi)modular lattices it cuts short the search
tree early.  The cutoff is simple to implement, but has a great impact
on the running time, which now scales roughly as~${2.8}^n$ for modular
lattices.  All vertically indecomposable modular $24$-lattices are
generated in about three minutes of processor time, compared to
(estimated) two years with Jipsen and Lawless's program.

With this faster algorithm modular lattices were counted up to size
$30$.  This gives an independent verification of Jipsen and Lawless's
numbers up to size~$23$ and a correction to their number for
size~$24$.  Semimodular lattices were counted up to $25$ elements,
graded lattices up to $21$, and geometric lattices up to $34$.  The
relevant entries in the Online Encyclopedia of Integer
Sequences~\cite{oeis} are A006981 (modular), A229202 (semimodular),
A278691 (graded), and A281574 (geometric).

The program code and the generated lattices (in compressed
\emph{digraph6} format~\cite{nauty}) are available for
download~\cite{eudat1,eudat2,eudat3}.

\section{Preliminaries}
\label{sec:preliminaries}
All lattices considered here are finite.  A lattice that has $n$
elements is called an $n$-\emph{lattice}, and its elements are
labeled with integers $i=1,2,\ldots,n$, with $1$ denoting the top
element.

The \emph{level} of an element, denoted by $\ell(i)$, is its longest
distance from the top, thus $\ell(1)=0$, coatoms have level 1 and so
on.  Without loss of generality we assume that element numbering is
consistent with levels, so that if $\ell(i)<\ell(j)$, then also $i<j$,
where $<$ denotes numerical order.  For a lattice~$L$, the set of
elements on level $k$ is denoted by~$L_k$.  The \emph{length} of a
lattice is the length of its longest chain, or equivalently, the level
of its bottom element.

We write $a \succ b$ if $a$ covers $b$.  The \emph{upper cover} of an
element $b$ is the set $\{a \; : \; a \succ b\}$, and its \emph{lower
  cover} is $\{a \; : \; b \succ a\}$.  The \emph{up-degree} and
\emph{down-degree} of an element are the sizes of its upper and lower
cover, respectively.

A lattice $L$ is \emph{vertically decomposable} if there is an element
distinct from top and bottom and comparable to every element of $L$.
Otherwise $L$ is a \emph{vertically indecomposable} lattice
(\emph{vi-lattice} for short).  For counting purposes, we only need to
generate the vi-lattices, since the total numbers can then be
calculated with a recursion formula~\cite{heitzig2002}.

\section{Algorithm for graded lattices}
\label{sec:basic}
A lattice is \emph{graded} if every maximal chain has the same length.
We begin by describing our basic algorithm that generates all
unlabeled, vertically indecomposable graded lattices of at~most $N$
elements.  By ``unlabeled'' we mean that it lists exactly one
representative of each isomorphism class, although for practical
purposes the lattices are represented with labeled elements.

The algorithm begins with lattices of length~2, and then recursively
adds new levels to create lattices up to the desired size~$N$.  (We
assume that lattices of lengths 0 and 1 are handled separately.)  The
initial lattices are $M_2,\ldots,M_{N-2}$, where $M_j$~denotes the
lattice that consists of the top, $j$ coatoms, and the bottom.  The
three-element lattice $M_1$ is omitted since it is not vertically
indecomposable.

The recursive step takes graded ``mother'' lattices of length $k$, and
creates graded ``daughter'' lattices of length $k+1$.  Let $L$ be a
mother lattice of length $k$ (so its atoms are at level $k-1$).
Daughter lattices are constructed by creating new elements, one at a
time, at level $k$.  Creating a new element involves specifying its
upper cover as a subset of $L_{k-1}$.  The possible upper covers are
considered in order of decreasing size: first we consider a new
element covered by all of $L_{k-1}$, then new elements covered by
$|L_{k-1}|-1$ elements chosen from $L_{k-1}$, and so on, finally down
  to elements covered by a single element of $L_{k-1}$.  Upper covers
  of the same size are considered in lexicographic order.  Whenever
  creating a new element, the algorithm checks that the proposed
  element does not violate the lattice conditions; see, \emph{e.g.},
  \cite{jipsen2015} for details.

As a new element is created, the resulting daughter lattice may not be
graded, since some elements of $L_{k-1}$ may not yet have been
included in any upper cover.  Such non-graded daughter lattices are
accepted as an intermediate step, but when level $L_k$ is completed,
we require that all elements of $L_{k-1}$ have been used in an upper
cover, ensuring that the lattice is graded.  Also, each level apart
from top and bottom is required to contain at least two elements in
order to restrict to vertically indecomposable lattices.

The basic method of ensuring nonisomorphism uses the graph-isomorphism
tool \texttt{nauty}~\cite{mckay2014,nauty}.  When listing daughters of
a given mother lattice, each daughter is converted to canonical
labeling and stored, along with three hash keys computed with the
\texttt{hashgraph} function provided by \texttt{nauty}.  A newly
created daughter lattice, also converted to canonical labeling, is
checked against this list and rejected if it is identical to any
previous daughter of the same mother.  Note that such a list only
needs to contain the accepted daughters of one mother, since daughters
of different (nonisomorphic) mother lattices are automatically
nonisomorphic.  That is, the algorithm does not need to keep
\emph{all} generated lattices in memory.

The method described above is in principle sufficient.  However, to
reduce work we employ a few shortcuts, depending on the structure of
the mother lattice.  Let $L$ be the mother lattice and $L_{k-1}$ its
atoms.  By examining the orbits and generators of its automorphism
group, as given by \texttt{nauty}, we classify $L$ into one of the
following types.

\textbf{Type 1, ``fixed''.}  Each atom is a singleton orbit, that is,
none of the automorphisms of $L$ move any atom.  In this case we need
not explicitly test the daughter lattices for isomorphism.  Different
daughter lattices have different collections of subsets of $L_{k-1}$
as the upper covers of elements on $L_k$.  Since all elements of
$L_{k-1}$ are fixed in all automorphisms, any two daughters are
nonisomorphic.

\begin{figure}[bt]
  \centering
  \begin{tikzpicture}
    \matrix(a)[matrix of math nodes, column sep=0.33cm, row sep=0.7cm]{
      L_0 &&&&&&&&& \node(1){1}; \\
      L_1 &&& \node(2){2}; &&&&&& \node(3){3}; && \node(4){4}; && \node(5){5}; \\
      L_2 & \node(6){6}; && \node(7){7}; && \node(8){8}; && \node(9){9}; && \node(10){10}; && \node(11){11}; && \node(12){12}; \\
      L_3 &&&&&&&&& \node(13){13}; \\
    };
    \draw (1)--(2);
    \draw (1)--(3);
    \draw (1)--(4);
    \draw (1)--(5);
    \draw (2)--(6);
    \draw (2)--(7);
    \draw (2)--(8);
    \draw (3)--(9);
    \draw (3)--(10);
    \draw (3)--(11);
    \draw (4)--(11);
    \draw (5)--(11);
    \draw (5)--(12);
    \draw (6)--(13);
    \draw (7)--(13);
    \draw (8)--(13);
    \draw (9)--(13);
    \draw (10)--(13);
    \draw (11)--(13);
    \draw (12)--(13);
    \draw[red,very thick] ($(6.north west)+(-0.1,0.1)$) rectangle ($(8.south east)+(0.1,-0.1)$);
    \draw[red,very thick] ($(9.north west)+(-0.1,0.1)$) rectangle ($(10.south east)+(0.1,-0.1)$);
    \draw[red,very thick] ($(11.north west)+(-0.1,0.1)$) rectangle ($(11.south east)+(0.1,-0.1)$);
    \draw[red,very thick] ($(12.north west)+(-0.1,0.1)$) rectangle ($(12.south east)+(0.1,-0.1)$);
  \end{tikzpicture}
  \caption{A mother lattice of the ``simple'' type, with rectangles
    indicating symmetric boxes.  For the first element on $L_3$ in a
    daughter lattice, the algorithm will, for example, consider the
    cover $6,7,9,12$ (consisting of prefixes of each box), but will
    not consider the cover $6,8,10,12$ which would be create an
    isomorphic copy.}
  \label{fig:simple}
\end{figure}
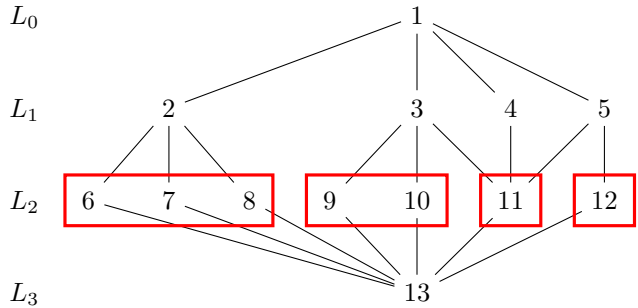
 
\textbf{Type 2, ``simple''.}  Some atoms are not singleton orbits, but
$L_{k-1}$ can be partitioned into subsets, \emph{symmetric boxes},
such that the elements in each box are in fully symmetric position.
To be more precise, $B \subseteq L_{k-1}$ is a symmetric box if, for
any permutation of $B$, there is an automorphism of~$L$ that moves $B$
by that permutation and leaves all other atoms fixed; and furthermore
$B$ is maximal in this respect.  If the mother lattice is of the
``simple'' type, then creating the first element $a$ of the next level
($L_k$) proceeds as follows.  Instead of considering all subsets $U
\subseteq L_{k-1}$ as possible upper covers of~$a$, we require that
for each symmetric box $B \subseteq L_{k-1}$, the intersection $U \cap
B$ is a prefix of $B$ in the numerical order of the elements.  An
empty prefix is allowed.  For example, if one of the symmetric boxes
is $B = \{6,7,8\}$, then $U \cap B$ may be $\varnothing$, $\{6\}$,
$\{6,7\}$, or $\{6,7,8\}$, but not $\{6,8\}$.  The same requirement is
held for each symmetric box, so the upper cover $U$ must be a union of
such numerical prefixes (some of which may be empty).  This is
illustrated in Fig.~\ref{fig:simple}.  This requirement drastically
cuts down the number of different upper covers that we need to
consider, especially if some of the symmetric boxes are large.  It
does not lead to missing any isomorphism classes, since for any
lattice $L'$ that does not fulfill this requirement, the algorithm
will visit a lattice that does, and is isomorphic to $L'$.  For
subsequent elements on each level we use the canonical labeling
method.

\textbf{Type 3, ``other''.}  In this case we employ no shortcuts but
simply use the canonical labeling method described above.

In practice the majority of mother lattices encountered fall into the
first two classes.  For example, among vertically indecomposable
graded $15$-lattices (there are $372\;838$ of them) the proportions of
the three types are 70.6\% fixed, 28.5\% simple, and 0.9\% other, so
in most cases some shortcuts apply.  However, it should be noted that
these shortcuts are quite simple, and a more carefully designed method
might reduce the amount of work even further.

Our basic isomorphism-avoidance method is similar to, but subtly
different from that used by Jipsen and Lawless~\cite{jipsen2015}.
Their approach is basically to create every daughter lattice, find a
canonical labeling with \texttt{nauty}, and then accept the daughter
if and only if it is the canonical daughter; this is checked by
inspecting its canonical labeling.  This ensures that from every
isomorphism class, exactly one lattice is accepted.  The benefit of
their approach is that the newly created lattice need not be compared
to any previously created lattices, so the previous lattices need not
be kept in memory (the so-called \emph{orderly} method).  However, for
this approach to work, one has to make sure that the canonical
daughter is indeed created.  For our approach it suffices that at
least \emph{one} representative (not necessarily the canonical one) of
each isomorphism class is created.  This allows some freedom in
designing shortcuts such as those described above.  If large numbers
of daughter lattices are not visited at all, the savings from this can
more than offset the extra work of memorizing the accepted lattices
and searching among them; the search is very fast anyway with the help
of hash tables.

On each level, some further optimizations are applied in the final
phase when creating elements of up-degree one.  They are not created
one by one; instead, for each element $a \in L_{k-1}$, we create some
\emph{number} $m(a)$ of elements on $L_k$, each of which is covered by
$a$ only.  We iterate over the possible choices of these integers
$m(a) \ge 0$, subject to the constraint that their sum does not cause
the lattice size to exceed the specified size~$N$, and further
requiring that $m(a) \ge 1$ for such $a$ whose lower covers are still
empty (otherwise the resulting lattice would not be graded).  For
details we refer to the program code~\cite{eudat1}.

\section{Algorithm for semimodular and modular lattices}
\label{sec:modular}

A finite lattice is \emph{semimodular} if, whenever $a \ne b$ and $a,b
\succ d$, there is an element $c$ such that $c \succ a,b$.  Dually, a
finite lattice is \emph{lower semimodular} if, whenever $c \succ a,b$
and $a \ne b$, there is an element $d$ such that $a,b \succ d$.
A~finite lattice that fulfills both conditions is \emph{modular}.
Note that the initial lattices $M_j$ in our recursive algorithm are
modular.  All semimodular and lower semimodular lattices are
graded~\cite[\S 3.3]{stanley1986}.

To generate semimodular lattices, we employ the graded lattice
algorithm from the previous section with two added conditions.  As
noted in the previous section, the algorithm creates new elements on
level $L_k$ by choosing, for each new element, an upper cover $U
\subseteq L_{k-1}$.  Here we additionally require that any two
elements $a,b \in U$ have a common covering element $c \succ a,b$.
The other condition is checked at the end of constructing a lattice:
we require that any two atoms have a common covering element.  (We
have to check this separately because the program does not explicitly
create the bottom element.)

To generate lower semimodular lattices, after the $k$th level of
elements is completed, we check that any pair of elements on $L_{k-1}$
that has a common covering element on $L_{k-2}$ has also been assigned
a common covered element on~$L_k$.

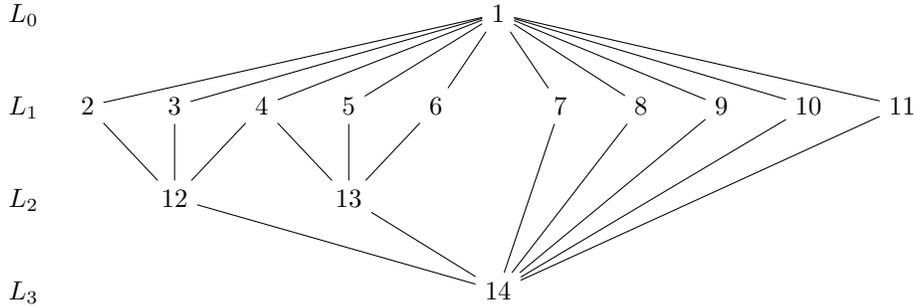
\begin{figure}[bt]
  \centering
  \begin{tikzpicture}
    \matrix(a)[matrix of math nodes, column sep=0.33cm, row sep=0.7cm]{
      L_0 &&&&&&&&&& \node(1){1}; \\
      L_1 &\node(2){2}; && \node(3){3}; && \node(4){4}; && \node(5){5}; && \node(6){6}; &&
      \node(7){7}; && \node(8){8}; && \node(9){9}; && \node(10){10}; && \node(11){11}; \\
      L_2 &&& \node(12){12}; &&&& \node(13){13}; \\
      L_3 &&&&&&&&&& \node(14){14}; \\
    };
    \draw (1)--(2);
    \draw (1)--(3);
    \draw (1)--(4);
    \draw (1)--(5);
    \draw (1)--(6);
    \draw (1)--(7);
    \draw (1)--(8);
    \draw (1)--(9);
    \draw (1)--(10);
    \draw (1)--(11);
    \draw (2)--(12);
    \draw (3)--(12);
    \draw (4)--(12);
    \draw (4)--(13);
    \draw (5)--(13);
    \draw (6)--(13);
    \draw (7)--(14);
    \draw (8)--(14);
    \draw (9)--(14);
    \draw (10)--(14);
    \draw (11)--(14);
    \draw (12)--(14);
    \draw (13)--(14);
  \end{tikzpicture}
  \caption{Illustration of the pair budget cutoff for lower
    semimodular lattices.}
  \label{fig:pairbudget}
\end{figure}

Again, the basic method described above is in principle sufficient to
generate the lower semimodular lattices, but a lot of work can be
avoided by an early cutoff that we will call the \emph{pair budget}.
We begin with an introductory example.  Consider the situation in
Fig.~\ref{fig:pairbudget}, where $L_1$ contains $10$ elements, and on
$L_2$ so far two elements have been created, both with up-degree~3.
Suppose further that we are listing lattices of at most $25$ elements.
On $L_1$ there are $\binom{10}{2}=45$ unordered pairs of distinct
elements.  For each such pair $a,b$, because $1 \succ a,b$, then for
lower semimodularity there must exist $d \in L_2$ such that
$a,b \succ d$.  Six pairs are already taken care of by the elements
labeled 12 and 13, so $45-6=39$ pairs are still wanting.  But recall
that new elements are added in order of decreasing up-degree.  Thus
any remaining element to be introduced on $L_2$ will have up-degree of
either~3 (in which case it is covered by three pairs on $L_1$), or~2
(covered by one pair), or~1 (covered by no pairs).  Since at most
$25-14=11$ more elements can be added on $L_2$, they will take care of
at most $11 \times 3 = 33$ pairs on $L_1$.  Clearly it is not possible
to extend this lattice into a lower semimodular one within the budget
of $25$ elements, and this branch of the search can be cut off
immediately.

In general the pair budget cutoff works as follows.  When beginning
level $L_k$, we first count the distinct pairs $a,b \in L_{k-1}$ such
that there exists $c \succ a,b$.  This is the number of pairs that
needs to be ``taken care of''.  Then, whenever a new element of
up-degree $r$ is created on $L_k$, we observe that it is covered by
$\binom{r}{2}$ pairs on $L_{k-1}$.  Any remaining element on $L_k$
will have up-degree of $r$ or smaller, and is thus covered by at most
$\binom{r}{2}$ pairs.  If, considering the maximum allowed size of a
lattice, the remaining elements cannot take care of enough pairs, this
branch is cut off.

The savings from the pair budget cutoff can be quite large.  Consider
again the situation in Fig.~\ref{fig:pairbudget}.  If the search were
not cut off, it would be possible to create daughter lattices where
each of the remaining $11$ elements chooses one of the remaining $39$
pairs as its upper cover, producing $\binom{39}{11} \approx 1.68
\times 10^{9}$ daughters.  More daughters would be created including
elements of up-degrees 3 and~1.  The actual number of daughters
visited would be somewhat smaller due to isomorphism.  But from the
simple counting argument we already know that \emph{none} of these
daughters can be lower semimodular.  Empirically, the total running
time for generating modular vi-lattices of 21 elements is cut more
than $40$-fold just by the pair budget cutoff, and the effect grows
with increasing lattice size.

An alternative way of generating semimodular lattices is to generate
lower semimodular lattices and then take their duals.  With the pair
budget method, this turned out to be much faster than generating
semimodular lattices directly, so this was the method we applied for
counting semimodular lattices.

In order to generate modular lattices, we simply use the algorithm
with both conditions (semimodularity and lower semimodularity).

\section{Algorithm for geometric lattices}
\label{sec:geometric}

A lattice is \emph{atomistic} if every element is a join of atoms, or
equivalently, if every element whose down-degree equals one is an
atom. (Stanley calls them \emph{atomic}, but we avoid this usage as
\emph{atomic} has other meanings.)  A finite lattice is
\emph{geometric} if it is semimodular and atomistic~\cite[\S
3.3]{stanley1986}.

There do not seem to be any previous computational approaches to
generating or counting geometric lattices, except that the present
author counted them up to size~15 (sequence A281574 in~\cite{oeis})
just by selecting geometric lattices from the lattice listings made
public by Malandro~\cite{malandro}.

Our algorithm actually generates the duals of geometric lattices, that
is, lower semimodular coatomistic lattices.  We use the algorithm for
lower semimodular lattices, with the extra condition that every
element below the coatom level ($L_1$) must have up-degree greater
than one.

\section{Partial verification}

Any attempt to establish mathematical truth by computation is prone to
many kinds of errors: hardware failures, human errors in operating the
computation, and errors in the algorithms and their software
implementation.  For example, Heitzig and Reinhold~\cite{heitzig2002}
observed wrong results in an earlier counting of unlabeled lattices;
and Brinkmann and McKay~\cite{brinkmann2002}, when counting posets of
up to 16 elements, experienced recurring hardware errors and had to
exclude some machines from their computations.

Short of formal verification of software and hardware, the reliability
of computational results rests on general qualities of the process,
such as simplicity, transparency and repeatability, and on various
kinds consistency checks.  In this work, several consistency checks
were employed to partially verify the results.  We can envisage three
types of errors in our lattice lists.  The lists might be incomplete;
they might contain objects that are not valid for the relevant lattice
family; or they might contain isomorphic duplicates of the same
lattice.

The first test, mainly against hardware errors, is that all lattice
lists were generated twice on separate computers of different models.
The resulting lattice lists (as text files) were verified to have the
same MD5 hash values.  Of course, any amount of repetition would not
help against logical errors in the program itself.

The second test looks for invalid objects in the listings.  The
lattices were checked with a separate program
(\texttt{latgrep}~\cite{eudat1}) to be of the relevant type (for
example, modular lattices).

The third test looks for isomorphic duplicates.  Each lattice list was
verified to be free of isomorphs by converting to a canonical labeling
with the \texttt{labelg} tool from the \texttt{nauty} package, and
then checking that all lines of the text file are differt.
A~weakness of this method is that it relies on the same
graph-isomorphism library (\texttt{nauty}) as the code that generated
the lattices.

The fourth test is between lattice families.  For sizes up to 25, we
verified that our lists of modular lattices are identical to what is
obtained by \emph{selecting} modular lattices from our lists of
semimodular lattices, with a separate filtering program
(\texttt{latgrep}).  A similar comparison was performed by selecting
geometric lattices from semimodular lattices (up to size 25), and
semimodular lattices from graded lattices (up to size 21).

The fifth test is by duality.  The families of modular and graded
lattices are closed with respect to duality.  For each lattice that we
listed for these two families, we checked that its dual (after
canonical labeling) also appears in the same listing.  Since the
generating algorithm builds the lattices in an inherently asymmetric
fashion from top to bottom, it seems plausible that errors in the
generatic logic would have been caught by failing the duality test.
However, this test would not detect missing self-dual lattices.

The sixth test is by comparison to previously published results.  The
counts match those computed by Jipsen and Lawless for modular lattices
up to 23, and for semimodular lattices up to 22
elements~\cite{jipsen2015}.  The numbers do \emph{not} match for
modular 24-lattices (we list exactly one more lattice).  Due to our
several consistency checks we are confident that the previous result
is in error.  Concerning actual lattice lists, Jipsen and Lawless's
program was rerun to generate modular and semimodular vi-lattices up
to 21 elements; after converting to canonical form with
\texttt{nauty}, the lists are identical to ours.  Unfortunately no
lattice listing from the previous result for modular 24-lattices was
available for comparison.

Another comparison to previous results concerns distributive lattices.
Since distributive lattices are modular, we can select distributive
vi-lattices of up to $30$ elements from our lists of modular
vi-lattices.  The counts thus obtained match those provided by Ern\'e
\emph{et al.}~\cite[Table 1]{erne2002}.

\section{Performance}

\begin{figure}[tb]
  \includegraphics[width=\textwidth]{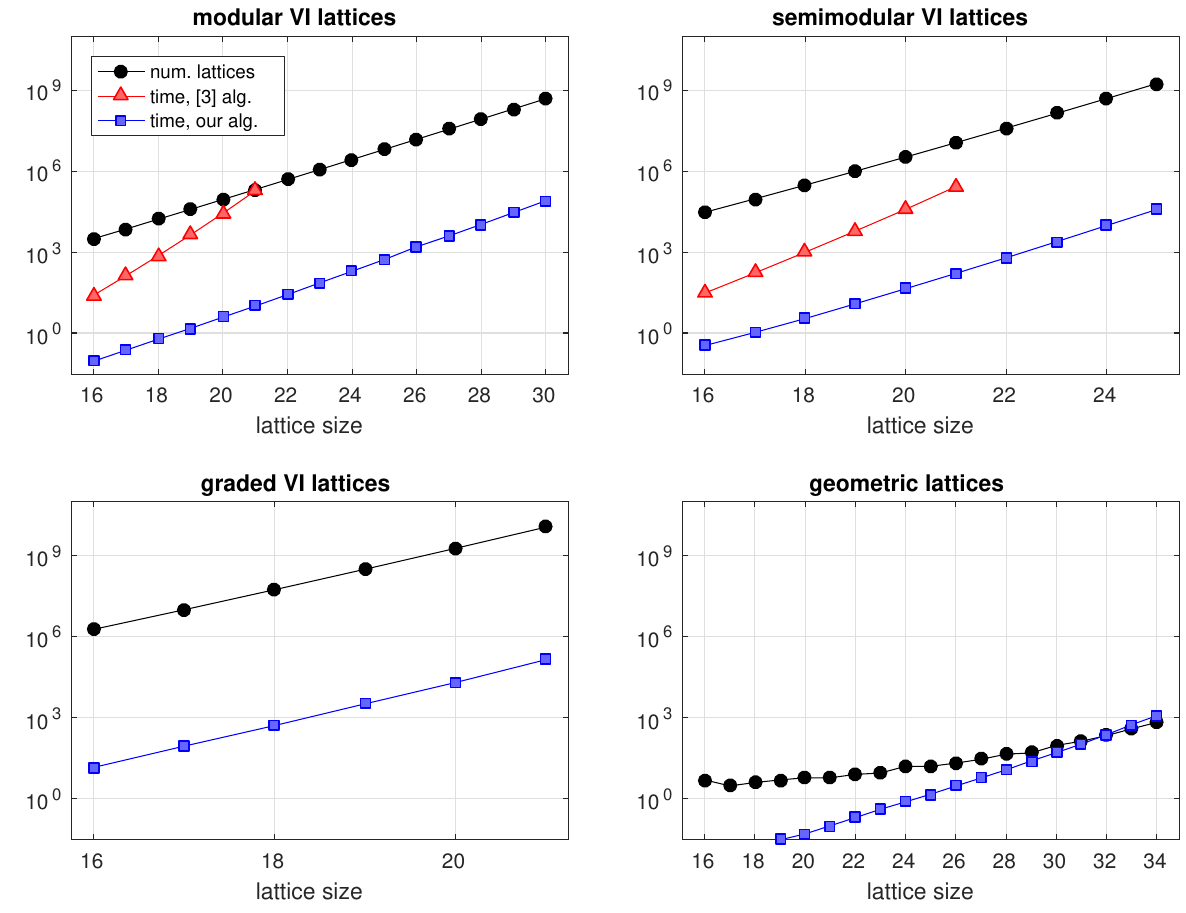}
  \captionsetup{width=.9\linewidth}
  \caption{Number of vertically indecomposable lattices by size, and
    running times of two algorithms (\cite{jipsen2015} and ours).  All
    times are in seconds on a single 2.6~GHz Intel Xeon \mbox{E5-2690}
    core.}
  \label{fig:timeplots}
\end{figure}

We do not have theoretical guarantees on the running time of our
algorithm, but some empirical observations can be made.
Fig.~\ref{fig:timeplots} illustrates the number of lattices and the
time spent by our algorithm.  Both exhibit somewhat similar scaling
with respect to lattice size, indicating that the algorithm is doing a
reasonable job in finding the relevant portions of the search space
(except for geometric lattices).  Let us inspect in more detail the
numerical growth ratios between consecutive lattice sizes.

\textbf{Modular vi-lattices.}  Between sizes $n=27$, $28$, $29$, $30$,
the number of vertically indecomposable modular lattices grows by
ratios $2.38$, $2.38$, $2.39$ (see Table~\ref{table:modsemi}),
suggesting a rather stable exponential growth.  For the same sizes our
running time grows by ratios $2.68$, $2.75$, $2.74$.  So empirically
the growth rate is slightly below ${2.4}^n$ for number of lattices and
${2.8}^n$ for running time.  This is not quite ideal: the gap of about
$0.4$ between the bases of the exponents raises the question whether
one could construct an \emph{output-sensitive} algorithm to generate
modular lattices, that is, one whose running time is linear in the
size of the output.

\textbf{Semimodular vi-lattices.}  Across sizes $n=22$, $23$, $24$,
$25$, the number of lattices grows roughly as~${3.6}^n$ and our
running time as~${4.0}^n$, again with a gap of~$0.4$ between the
bases.

\textbf{Graded vi-lattices.}  Across sizes $n=18$, $19$, $20$, $21$,
the number of lattices grows by ratios $5.74$, $5.93$, $6.13$,
exhibiting faster than exponential growth.  The running time grows by
ratios $6.38$, $6.16$ and $7.03$, again somewhat faster than the
number of lattices.

\textbf{Geometric lattices.}  The number of geometric lattices grows
so slowly that no asymptotic form is discernible from the available
numbers.  The running time grows much faster than the number of
lattices.

To compare the performance, we reran Jipsen and Lawless's
program~\cite{jipsen2015} for modular and semimodular vi-lattices,
both up to size 21.  These running times are shown in
Fig.~\ref{fig:timeplots} with red triangles.  Between sizes $n=18$,
$19$, $20$, $21$, the running time for modular vi-lattices grew by
ratios $5.90$, $6.27$, $6.81$, with $n=21$ taking $2.1$ days of
processor time.  From this we estimate that $n=24$ would have required
about two years of processor time (about $300\;000$ times more than
with our algorithm, which completed in 194 seconds).

We conclude this section with some remarks on the speed of our basic
algorithm for generating graded lattices.  At size $n=21$ it outputs
about $82\;000$ graded lattices per second, or one lattice in $12$
microseconds.  On the 2.6~GHz processor that was used, this amounts to
$32\;000$ clock cycles per lattice.  Gebhardt and
Tawn~\cite{gebhardt2016} count general (not graded) lattices
considerably faster, at $2\;200$ clock cycles per lattice.  They
handle isomorphisms with a sophisticated method specially tailored to
lattices.  In contrast, our algorithm handles only the simplest cases
directly, and in more complicated cases resorts to using the
general-purpose graph isomorphism tool \texttt{nauty} as a ``black
box''.  Indeed, during the search for graded 21-lattices, our
algorithm performs about $2.7 \times 10^{10}$ calls to \texttt{nauty}.
Taking $3.4$ microseconds per call on average, together they account
for 65\% of the total running time.  While Gebhardt and Tawn
considered general lattices only, it might be useful to combine their
method for isomorphisms with our early checks for (semi)modularity.

\section{Results}

The lattice listings are available for
download~\cite{eudat1,eudat2,eudat3}.  The numbers of lattices are
shown in Tables \ref{table:modsemi} and~\ref{table:grageo}.
Vertically indecomposable modular, semimodular and graded lattices
were directly counted by the program; numbers that include
decomposable lattices were then calculated with the recursion
formula~\cite{heitzig2002}
\[
u(n) = \sum_{k=2}^n u_\text{vi}(k) \; u(n\!-\!k\!+\!1),
\qquad\text{for $n\ge 2$},
\]
where $u(n)$ counts all unlabeled lattices of size~$n$ in the relevant
family, and $u_\text{vi}(n)$ counts vi-lattices only.  For geometric
lattices the recursion formula does not apply as they are necessarily
vertically indecomposable.

\begin{figure}[tb]
  \begin{minipage}[t]{0.48\linewidth}
    \includegraphics[width=\textwidth]{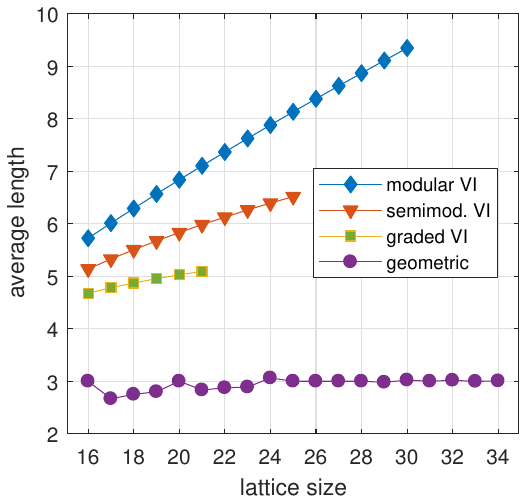}
    \captionsetup{width=.9\linewidth}
    \caption{Average lattice length as a function of lattice size.}
    \label{fig:avelen}
  \end{minipage}
  \hfill
  \begin{minipage}[t]{0.48\linewidth}
    \includegraphics[width=\textwidth]{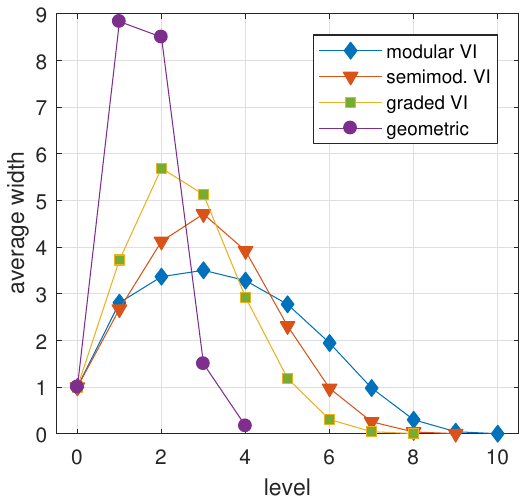}
    \captionsetup{width=.9\linewidth}
    \caption{Average widths (numbers of elements) of levels in
      lattices of size~21.}
    \label{fig:avewid}
  \end{minipage}
\end{figure}

The numbers seem to suggest exponential growth of modular and
semimodular lattices; indeed, Jipsen and Lawless~\cite{jipsen2015}
have proven a lower bound of the form $\mathrm\Omega(2^n)$ for the
number of modular lattices.  This raises the question of finding an
upper bound $O(c^n)$ with some constant~$c$.  To the best of our
knowledge, no such upper bound is known for modular or semimodular
lattices.  In contrast, it is known that the number of graded lattices
grows faster than exponentially in~$n$~\cite{kleitman1980}.

Apart from total numbers, one may compute various statistics from the
actual lattice listings.  As an example, from Figs.~\ref{fig:avelen}
and~\ref{fig:avewid} we observe that typical lattices in these four
families have quite different length and width characteristics:
modular lattices are long and narrow while geometric lattices are
short and wide.  Semimodular and graded lattices are in between.
Empirical understanding of the typical lattice shape may be useful,
for example, when formulating hypotheses about asymptotics, and in
algorithm design.

\begin{table}[p]
\centering
\begin{tabular}{r|crrc|crr}
  $n$ && modular vi & modular &&& semimodular vi & semimodular \\
  \hline
  1	&& 1	& 1	&&& 1	& 1 \\
  2	&& 1	& 1	&&& 1	& 1 \\
  3	&& 0	& 1	&&& 0	& 1 \\
  4	&& 1	& 2	&&& 1	& 2 \\
  5	&& 1	& 4	&&& 1	& 4 \\
  6	&& 2	& 8	&&& 2	& 8 \\
  7	&& 3	& 16	&&& 4	& 17 \\
  8	&& 7	& 34	&&& 9	& 38 \\
  9	&& 12	& 72	&&& 21	& 88 \\
  10	&& 28	& 157	&&& 53	& 212 \\
  11	&& 54	& 343	&&& 139	& 530 \\
  12	&& 127	& 766	&&& 384	& 1\;376 \\
  13	&& 266	& 1\;718	&&& 1\;088	& 3\;693 \\
  14	&& 614	& 3\;899	&&& 3\;186	& 10\;232 \\
  15	&& 1\;356	& 8\;898	&&& 9\;596	& 29\;231 \\
  16	&& 3\;134	& 20\;475	&&& 29\;601	& 85\;906 \\
  17	&& 7\;091	& 47\;321	&&& 93\;462	& 259\;291 \\
  18	&& 16\;482	& 110\;024	&&& 301\;265	& 802\;308 \\
  19	&& 37\;929	& 256\;791	&&& 990\;083	& 2\;540\;635 \\
  20	&& 88\;622	& 601\;991	&&& 3\;312\;563	& 8\;220\;218 \\
  21	&& 206\;295	& 1\;415\;768	&&& 11\;270\;507	& 27\;134\;483 \\
  22	&& 484\;445	& 3\;340\;847	&&& 38\;955\;164	& 91\;258\;141 \\
  23	&& 1\;136\;897	& 7\;904\;700	&&& \cellcolor{blue!15}136\;660\;780	& \cellcolor{blue!15}312\;324\;027 \\
  24	&& \cellcolor{blue!15}2\;682\;451\rlap{$^*$}	& \cellcolor{blue!15}18\;752\;943\rlap{$^*$}	&&& \cellcolor{blue!15}486\;223\;384	& \cellcolor{blue!15}1\;086\;545\;705 \\
  25	&& \cellcolor{blue!15}6\;333\;249	& \cellcolor{blue!15}44\;588\;803	&&& \cellcolor{blue!15}1\;753\;185\;150	& \cellcolor{blue!15}3\;838\;581\;926 \\
  26	&& \cellcolor{blue!15}15\;005\;945	& \cellcolor{blue!15}106\;247\;120 & \\
  27	&& \cellcolor{blue!15}35\;595\;805	& \cellcolor{blue!15}253\;644\;319 & \\
  28	&& \cellcolor{blue!15}84\;649\;515	& \cellcolor{blue!15}606\;603\;025 & \\
  29	&& \cellcolor{blue!15}201\;560\;350	& \cellcolor{blue!15}1\;453\;029\;516 & \\
  30	&& \cellcolor{blue!15}480\;845\;007	& \cellcolor{blue!15}3\;485\;707\;007 & \\
  \hline
\end{tabular}
\caption{Numbers of unlabeled modular and semimodular lattices by size
  (vi~=~vertically indecomposable).  New numbers are highlighted;
  corrections to previous numbers marked with $^*$.}
\label{table:modsemi}
\end{table}

\begin{table}[p]
\centering
\begin{tabular}{r|crrc|cr}
  $n$  && graded vi & graded &&& geometric \\
  \hline
  1	&& 1	& 1	&&& 1 \\
  2	&& 1	& 1	&&& 1 \\
  3	&& 0	& 1	&&& 0 \\
  4	&& 1	& 2	&&& 1 \\
  5	&& 1	& 4	&&& 1 \\
  6	&& 3	& 9	&&& 1 \\
  7	&& 7	& 22	&&& 1 \\
  8	&& 22	& 60	&&& 2 \\
  9	&& 68	& 176	&&& 1 \\
  10	&& 242	& 565	&&& 2 \\
  11	&& 924	& 1\;980	&&& 1 \\
  12	&& 3\;793	& 7\;528	&&& 3 \\
  13	&& 16\;569	& 30\;843	&&& 2 \\
  14	&& 76\;638	& 135\;248	&&& 2 \\
  15	&& 372\;838	& 630\;004	&&& 3 \\
  16	&& \cellcolor{blue!15}1\;900\;132	& \cellcolor{blue!15}3\;097\;780	&&& \cellcolor{blue!15}5 \\
  17	&& \cellcolor{blue!15}10\;105\;175	& \cellcolor{blue!15}15\;991\;395	&&& \cellcolor{blue!15}3 \\
  18	&& \cellcolor{blue!15}55\;895\;571	& \cellcolor{blue!15}86\;267\;557	&&& \cellcolor{blue!15}4 \\
  19	&& \cellcolor{blue!15}320\;655\;822	& \cellcolor{blue!15}484\;446\;620	&&& \cellcolor{blue!15}5 \\
  20	&& \cellcolor{blue!15}1\;903\;047\;753	& \cellcolor{blue!15}2\;822\;677\;523	&&& \cellcolor{blue!15}6 \\
  21	&& \cellcolor{blue!15}11\;658\;925\;558	& \cellcolor{blue!15}17\;017\;165\;987	&&& \cellcolor{blue!15}6 \\
  22	&& 	& 	&&& \cellcolor{blue!15}8 \\
  23	&& 	& 	&&& \cellcolor{blue!15}9 \\
  24	&& 	& 	&&& \cellcolor{blue!15}16 \\
  25	&& 	& 	&&& \cellcolor{blue!15}16 \\
  26	&& 	& 	&&& \cellcolor{blue!15}21 \\
  27	&& 	& 	&&& \cellcolor{blue!15}29 \\
  28	&& 	& 	&&& \cellcolor{blue!15}45 \\
  29	&& 	& 	&&& \cellcolor{blue!15}50 \\
  30	&& 	& 	&&& \cellcolor{blue!15}95 \\
  31	&& 	& 	&&& \cellcolor{blue!15}136 \\
  32	&& 	& 	&&& \cellcolor{blue!15}220 \\
  33	&& 	& 	&&& \cellcolor{blue!15}392 \\
  34	&& 	& 	&&& \cellcolor{blue!15}680 \\
\end{tabular}
\caption{Numbers of unlabeled graded and geometric lattices by size
  (vi~=~vertically indecomposable).  New numbers are highlighted.}
\label{table:grageo}
\end{table}

\section{Conclusion}

Let us conclude with some thoughts of future research.  Based on the
empirical running times and memory usage, it seems quite feasible to
count modular lattices a little further with our program as it is:
modular $32$-lattices might be counted in about one week of processor
time, requiring a few gigabytes of memory to keep the daughter lattice
lists for isomorph rejection.  The computation can be parallelized by
running a separate job for each initial lattice $M_j$.  Extending the
counts of semimodular lattices seems also feasible.

However, it might be more interesting to pursue other methods of
counting modular lattices, possibly without explicitly generating
them.  Jipsen and Lawless~\cite{jipsen2015} discuss some prospects of
using alternative representations of modular lattices to count them.
Another prospect comes from our observation that modular lattices tend
to be long and narrow: perhaps a large portion of modular lattices
could be counted implicitly, by considering some kind of vertical
compositions of smaller modular lattices, without explicitly listing
the compositions.

For counting graded lattices, the current program does not seem well
suited to larger sizes.  Compared to modular lattices, graded lattices
tend to be shorter and wider, which becomes a problem in our method of
isomorph rejection: a single mother lattice may have a large number of
daughter lattices, so the memory required for storing the daughters
becomes excessive.  It would seem necessary to improve the isomorph
rejection method so~as to use less memory.

Geometric lattices share the problem of being short and wide.  Another
problem is that our current method generates a large number of lower
semimodular proposals, then rejects the vast majority of them for not
being co-atomistic.  This seems very inefficient, and better methods
would be desirable, for example, ones that would reject the proposed
lattices earlier.

\section*{Acknowledgements}
The author wants to thank Nathan Lawless for providing the program
code described in~\cite{jipsen2015}, and the anonymous referee for
several valuable remarks.

The research that led to these results has received funding from the
European Research Council under the European Union's Seventh Framework
Programme (FP/2007-2013) / ERC Grant Agreement 338077 ``Theory and
Practice of Advanced Search and Enumeration.''

Computational resources were provided by CSC -- IT Center for Science,
Finland, and the Aalto Science-IT project.


\bibliographystyle{plain}
\bibliography{refs}

\end{document}